РОССИЙСКАЯ АКАДЕМИЯ НАУК
ФЕДЕРАЛЬНОЕ ГОСУДАРСТВЕННОЕ БЮДЖЕТНОЕ
УЧРЕЖДЕНИЕ НАУКИ
**ВЫЧИСЛИТЕЛЬНЫЙ ЦЕНТР**
**им. А.А.Дородницына**

**М.П. ВАЩЕНКО, А.В. ГАСНИКОВ, Е.Г. МОЛЧАНОВ,**
**Л.Я. ПОСПЕЛОВА, А.А. ШАНАНИН**

# ВЫЧИСЛИМЫЕ МОДЕЛИ И ЧИСЛЕННЫЕ МЕТОДЫ ДЛЯ АНАЛИЗА ТАРИФНОЙ ПОЛИТИКИ ЖЕЛЕЗНОДОРОЖНЫХ ГРУЗОПЕРЕВОЗОК

Модель конкурентного равновесия модифицирована для анализа проблем формирования тарифной и инвестиционной политики управления железнодорожными грузоперевозками в современных российских условиях. Получен вариационный принцип в форме пары взаимно двойственных задач выпуклого программирования, из решения которых находится конкурентное равновесие. Это даёт возможность анализа актуальных проблем системы железнодорожных перевозок. Предложен подход к анализу перекрёстного субсидирования перевозки грузов с различной доходностью.

**Работа выполнена при финансовой поддержке РФФИ (проект 13-07-13110 офи_м_РЖД).**

Ключевые слова: общественные блага, модель грузоперевозок, коммуникационные ограничения, теорема Фенхеля, преобразование Лежандра-Юнга-Фенхеля, конкурентное равновесие.

Рецензенты: А.В. Арутюнов,
С.П. Тарасов

Научное издание

# 1. Введение

Результаты реформы ОАО «РЖД», которая началась в 2001 г., вызвали неоднозначную реакцию. С одной стороны, привлечение на рынок грузовых перевозок частных операторов существенно увеличило инвестиционную привлекательность отрасли. Удалось решить вопрос острого дефицита парка вагонов, возникший в начале 2000-х годов на фоне восстановления промышленности после кризиса 1998 г. С другой стороны, обнаружились негативные эффекты рыночных механизмов: резко выросли тарифы на грузовые железнодорожные перевозки, вырос порожний пробег и простой вагонов, что увеличило нагрузку на инфраструктуру.

Для сегментированной экономики России, в которой наряду с эффективными и высокорентабельными нефтегазовыми и частично горнорудными компаниями присутствуют депрессивные низко прибыльные предприятия обрабатывающего сектора, неограниченная либерализация, как правило, приводит к росту тарифов, инфляции потребительских цен в целом и неизбежным дополнительным социальным проблемам. Поэтому государственная политика часто в целях борьбы с инфляцией включает ограничение роста тарифов на услуги естественных монополий.

Однако ограничение тарифов негативно сказывается на развитии естественных монополий. Поэтому вопрос формирования тарифов на грузоперевозки крайне актуален для современной России не только с точки зрения самой железнодорожной инфраструктуры (увеличение-уменьшение трафика), но и с точки зрения среднесрочных экономико-социальных последствий.

Вопросам анализа воспроизводства и распределения товаров и услуг, предоставляемых государством и естественными монополиями, посвящена разработанная около ста лет назад концепция равновесия Э. Линдаля (1919) (см., например, [1-4]) в модели с общественными благами. Равновесие по Линдалю предполагает, что вложения экономического агента в финансирование общественного блага определяется предельной полезностью для него общественного блага. Под общественными благами в этой модели понимаются товары и услуги, которыми могут пользоваться одновременно разные потребители. Примерами отраслей, производящих такие товары и услуги, являются инфраструктурные отрасли (естественные монополии). Основным результатом исследования таких моделей является теорема о механизме конкурентного равновесия, который должен допускать дифференциацию цен (перекрестное субсидирование), т.е. потребители различных типов за предоставленное общественное благо должны платить по разным тарифам. Каким должен быть компромисс между интересами различных экономических агентов? Для анализа этого вопроса можно предложить разные подходы. Отметим, например, цикл работ по исследованию игр с иерархическим вектором интересов Гермейера-Вателя [5-9].

Известный парадокс Эджворта (см. [10]) показывает, что возможны неожиданные реакции при попытках регулирования тарифов на общественные блага.

Рассмотрим, например, грузоперевозчика, в собственности которого находятся вагоны и который устанавливает тарифы на перевозку вагонов двух типов между двумя пунктами. Предположим, что качество услуг и доходность перевозки грузов в вагонах разных типов различна. Она определяется тарифами $c_1$ на перевозку вагонов первого типа, $c_2$ на перевозку вагонов второго типа и оплатой услуг железнодорожной компании $t = 10$ (для вагонов обоих типов).

Предположим, что к грузоперевозчику обращаются клиенты четырёх типов. Клиенты первого типа обратятся с предложением перевезти 450 вагонов первого типа, если $c_1 - c_2 < 5, c_1 \le 22$, или перевезти груз вагонами второго типа, если $c_1 - c_2 \ge 5, c_2 \le 17$, при невыполнении этих условий они воспользуются альтернативными способами перевозки, например автотранспортом. Клиенты второго типа обратятся с предложением перевезти 40 вагонов первого типа, если $c_1 \le 21$, при невыполнении этого условия они воспользуются альтернативными способами перевозки, например автотранспортом.

Клиенты третьего типа обратятся с предложением перевезти 900 вагонов второго типа, если $c_2 \le 18$, при невыполнении этого условия они воспользуются альтернативными способами перевозки, например автотранспортом.

Клиенты четвертого типа обратятся с предложением перевезти 200 вагонов второго типа, если $c_2 \le 16$, при невыполнении этого условия они воспользуются альтернативными способами перевозки, например автотранспортом.

Нетрудно вычислить, что грузоперевозчик, который максимизирует свою прибыль, установит тарифы $c_1 = 22, c_2 = 18$.

Как изменятся тарифы, если железнодорожная компания повысит оплату своих услуг по перевозке вагонов первого типа до $18$, оставив оплату для вагонов второго типа прежней и равной 10?

Можно было бы ожидать повышения тарифов $c_1, c_2$. Однако несложные расчеты показывают, что максимальную прибыль перевозчик получит, снижая тарифы до $c_1 = 21, c_2 = 16$.

Пример Эджворта показывает, что подход к анализу стратегии отрасли, предоставляющей услуги по железнодорожным грузоперевозкам, должен основываться на вычислимых моделях взаимодействия экономических агентов. Удачным примером такого подхода может служить знаменитая работа Л.В.Канторовича и М.К.Гавурина [11], которая послужила исходной посылкой для наших исследований, адаптирующих классические математические модели экономических систем для анализа актуальных вопросов управления железнодорожными грузоперевозками в современных российских условиях.

Государственная политика сдерживания потребительской инфляции за счёт ограничения тарифов на услуги естественных монополий как следствие порождает дефицит финансовых ресурсов на обновление и модернизацию их основных фондов. В качестве возможного подхода к решению проблемы обычно предлагается выделить некоторый сегмент деятельности естественной монополии и привлечь в него частных инвесторов. В системе железнодорожных перевозок недостаток финансовых ресурсов привёл к необходимости решать проблему обновления вагонного парка.

Эта проблема была решена в процессе реформирования ОАО «РЖД» с привлечением частных инвесторов и допуском на рынок грузоперевозок коммерческих компаний, являющихся владельцами вагонного парка и посредниками между ОАО «РЖД» и потребителями услуг по грузоперевозке.

Появление таких посредников могло привести к росту тарифов на грузоперевозки, сокращению железнодорожного грузооборота и специализации на перевозке высокодоходных грузов. Поэтому при реформировании ОАО «РЖД» были выделены крупные компании ОАО «Первая грузовая компания», ОАО «Федеральная грузовая компания» и ОАО «ТрансКонтейнер». В управлении этими компаниями сохраняется «государственное влияние», структура их совокупного грузооборота в целом сохранила дореформенные пропорции, и они осуществляют более 80% грузоперевозок. Это позволило, обновив парк грузовых вагонов и устранив их дефицит, не допустить неконтролируемого роста тарифов и сокращения грузооборота. В частности, за счёт перекрёстного субсидирования сохраняется перевозка грузов с низкой доходностью, что важно для сохранения целостности страны. Однако проявились и негативные последствия. Небольшие коммерческие компании повысили тарифы, предложив лучшее качество услуг (более короткие сроки перевозок), и стали специализироваться на перевозке высокодоходных грузок. Такая специализация вызвала увеличение перегона порожняка, увеличила нагрузку на транспортную инфраструктуру в целом, что порождает увеличение сроков обслуживания. Кроме того, увеличение тарифов оказывает влияние на спрос на услуги по железнодорожным грузоперевозкам. В частности, стали распространяться перевозки высокодоходных грузов автотранспортом. Все эти последствия являются упущенной прибылью для системы железнодорожных грузоперевозок в целом и нуждаются в анализе, учитывающем косвенные последствия принимаемых решений.

Естественным языком для такого анализа являются математические модели, учитывающие специфику экономических отношений, сложившихся в отрасли. В этих моделях желательно учесть в явном виде влияние тарифов на спрос на услуги по железнодорожным грузоперевозкам и возможности извлечения сверхприбыли.

## 2. Конкурентное равновесие в модели железнодорожных грузоперевозок

Рассмотрим модель, в которой выделены $m$ видов товаров[1], пункты потребления товаров $I$, пункты производства $J$. Обозначим через $X_i^k, k \in K, K = \{1,...,m\}$ количество $k$ -го товара, поступившее в пункт потребления $i \in I$, а через $Y_j^k$ - количество $k$-го товара, вывезенное из пункта производства $j \in J$. Будем считать, что от поступивших товаров в пункте потребления $i \in I$ получается доход[2] $F_i\left(X_i^1,...,X_i^m\right)$, а себестоимость производства товаров в пункте производства $j \in J$ в количестве $\left(Y_j^1,...,Y_j^k\right)$ равна $G_j\left(Y_j^1,...,Y_j^m\right)$. Предположим, что эти функции удовлетворяют неоклассическим требованиям, т.е. непрерывны, монотонно не убывают по своим аргументам, функции $F_i\left(X_i^1,...,X_i^m\right)$ вогнуты, а функции $G_j\left(Y_j^1,...,Y_j^m\right)$ выпуклы. Обозначим через $c_{ij}^k$ затраты в денежном выражении на перевозку единицы $k$-го товара из пункта производства $j$ в пункт потребления $i$, через $z_{ij}^k$ - объем перевозок $k$-го товара из $j$-го пункта производства в $i$-й пункт потребления, через $p_i^k$ - цену $k$-го товара в $i$-м пункте потребления, через $\hat{p}_j^k$ - цену $k$-го товара в $j$-м пункте производства.

[1] При идентификации модели железнодорожных грузоперевозок, выделяемые товары могут классифицироваться по типам железнодорожных вагонов, используемых для их перевозки.

Рассмотрим сначала функционирование отрасли в условиях совершенной конкуренции, когда экономические агенты не могут влиять своим поведением на цены и максимизируют свою прибыль при заданных ценах.

**Определение.** *Будем говорить, что набор неотрицательных чисел*

$$\left\{X_i^k, Y_j^k, z_{ij}^k, p_i^k, \hat{p}_j^k \middle| i \in I, j \in J, k \in K\right\}$$

*является конкурентным равновесием в модели железнодорожных грузоперевозок, если*

*1) для любого* $i \in I$

---

[2] Функцию $F_i\left(X_i^1,...,X_i^m\right)$ можно интерпретировать как производственную функцию в $i$-м пункте потребления.

$$\left(X_i^1,...,X_i^m\right)\in Arg\max\left\{F_i\left(x_i^1,...,x_i^m\right)-\sum_{i=1}^{k}p_i^k x_i^k\,\middle|\left(x_i^1,...,x_i^m\right)\geq 0\right\},$$

*2) для любого* $j\in J$

$$\left(Y_j^1,...,Y_j^m\right)\in Arg\max\left\{\sum_{k=1}^{m}\hat{p}_j^k y_j^k-G_j\left(y_j^1,...,y_j^m\right)\middle|\left(y_j^1,...,y_j^m\right)\geq 0\right\},$$

*3) для любых* $i\in I,\ k\in K\ \ X_i^k\leq\sum_{j\in J}z_{ij}^k,\ p_i^k\left(X_i^k-\sum_{j\in J}z_{ij}^k\right)=0,$

*4) для любых* $j\in J,\ k\in K\ \ Y_j^k\geq\sum_{i\in I}z_{ij}^k,\ \hat{p}_j^k\left(Y_j^k-\sum_{i\in I}z_{ij}^k\right)=0,$

*5) для любых* $i\in I,\ j\in J$,

$$p_i^k\leq\hat{p}_j^k+c_{ij}^k,\ \ z_{ij}^k\geq 0,\ \ z_{ij}^k\left(\hat{p}_j^k+c_{ij}^k-p_i^k\right)=0.$$

**Замечание.** Предлагаемое понятие конкурентного равновесия в модели железнодорожных грузоперевозок соответствует вальрасовской концепции совершенного конкурентного равновесия, в которой субъекты экономики максимизируют свои прибыли, считая, что своими решениями они не влияют на сложившиеся цены на товары, а цены устанавливаются так, чтобы спрос на каждый вид товаров соответствовал его предложению. В отличие от традиционных постановок задач о грузоперевозках в модели учитывается эластичность спроса и предложения по ценам. Вычисление конкурентного равновесия может быть сведено к решению следующего вариационного неравенства (стандартной задачи дополнительности). Рассмотрим множество $R_+^N=\left\{\left(z_{ij}^k\,\middle|k\in K,i\in I,j\in J\right)\geq 0\right\}$. Определим многозначное отображение

$$P\left(z_{ij}^k \middle| k\in K, i\in I, j\in J\right)=$$

$$\left\{\left(\hat{p}_j^k - p_i^k + c_{ij}^k \middle| k\in K, i\in I, j\in J\right)\right|$$

$$\left.\left(\hat{p}_j^1,...,\hat{p}_j^m\right)\in\partial G\left(\sum_{i\in I} z_{ij}^k\right)-\partial F\left(\sum_{j} z_{ij}^k\right)+c_{ij}^k\right\}$$

из $R_+^N=\left\{\left(z_{ij}^k \middle| k\in K, i\in I, j\in J\right)\geq 0\right\}$ в

$R^N=\left\{\left(P_{ij}^k \middle| k\in K, i\in I, j\in J\right)\right\}$.

Здесь $\partial F_i(\bullet)$ - супердифференциал вогнутой функции $F_i(\bullet)$, $\partial G_i(\bullet)$ - субдифференциал выпуклой функции $G_i(\bullet)$.

Напомним, что решением стандартной задачи дополнительности $\left(R_+^N, P\left(z_{ij}^k \middle| k\in K, i\in I, j\in J\right)\right)$ называется точка $\left(\tilde{z}_{ij}^k \middle| k\in K, i\in I, j\in J\right)\in R_+^N$, для которой найдётся точка $\left(P_{ij}^k \middle| k\in K, i\in I, j\in J\right)\in P\left(z_{ij}^k \middle| k\in K, i\in I, j\in J\right)$, такая, что

$$P_{ij}^k\geq 0, P_{ij}^k\tilde{z}_{ij}^k=0 \ \left(k\in K, i\in I, j\in J\right).$$

Для эффективного анализа и численного решения вариационного неравенства $\left(R_+^N, P\left(z_{ij}^k \middle| k\in K, i\in I, j\in J\right)\right)$ построим пару взаимно двойственных задач выпуклой оптимизации, по решениям которых строится конкурентное равновесие.

Рассмотрим задачу о максимизации прибыли экономической системы с учетом затрат на грузоперевозки:

$$\sum_{i\in I}F_i\left(X_i^1,...,X_i^m\right)-\sum_{j\in J}G_j\left(Y_j^1,...,Y_j^m\right)-\sum_{\substack{i\in I,\\ j\in J,\\ k\in K}}c_{ij}^k z_{ij}^k\to\max \quad (2.1)$$

$$X_i^k \le \sum_{j\in J} z_{ij}^k \quad (i\in I, k\in K), \tag{2.2}$$

$$Y_j^k \ge \sum_{i\in I} z_{ij}^k \quad (j\in J, k\in K), \tag{2.3}$$

$$z_{ij}^k \ge 0 \quad (i\in I, j\in J, k\in K). \tag{2.4}$$

Функция прибыли $i$ -го потребителя равна

$$\Pi_i\left(p_i^1,...,p_i^m\right) = \sup\left\{F_i\left(X_i^1,...,X_i^m\right) - \sum_{k\in K} p_i^k X_i^k \,\middle|\, \left(X_i^1,...,X_i^m\right)\ge 0\right\}.$$

Функция прибыли $j$ -го производителя равна

$$\pi_j\left(\hat{p}_j^1,...,\hat{p}_j^m\right) = \sup\left\{\sum_{k\in K} \hat{p}_j^k Y_j^k - G_j\left(Y_j^1,...,Y_j^m\right) \,\middle|\, \left(Y_j^1,...,Y_j^m\right)\ge 0\right\}.$$

Нетрудно видеть, что функция $\Pi_i\left(p_i^1,...,p_i^m\right)$ является непрерывной, выпуклой, монотонно невозрастающей функцией по переменным $\left(p_i^1,...,p_i^m\right)$ на множестве $R_+^m$, а функция $\pi_j\left(\hat{p}_j^1,...,\hat{p}_j^m\right)$ является непрерывной, выпуклой, монотонно неубывающей функцией по переменным $\left(\hat{p}_i^1,...,\hat{p}_i^m\right)$ на множестве $R_+^m$.

По теореме Фенхеля-Моро [11] имеем

$$F_i\left(X_i^1,...,X_i^m\right) = \inf\left\{\Pi_i\left(p_i^1,...,p_i^m\right) + \sum_{k\in K} p_i^k X_i^k \,\middle|\, \left(p_i^1,...,p_i^m\right)\ge 0\right\}, \tag{2.5}$$

$$G_j\left(Y_j^1,...,Y_j^m\right) = \sup\left\{\sum_{k\in K} \hat{p}_j^k Y_j^k - \pi_j\left(\hat{p}_j^1,...,\hat{p}_j^m\right) \,\middle|\, \left(\hat{p}_j^1,...,\hat{p}_j^m\right)\ge 0\right\}. \tag{2.6}$$

Двойственной по Фенхелю [11, с. 46-47] экстремальной задачей к задаче (2.1)-(2.4) является задача выпуклого программирования:

$$\sum_{i\in I}\Pi_i\left(p_i^1,...,p_i^m\right)+\sum_{j\in J}\pi_j\left(\hat{p}_j^1,...,\hat{p}_j^m\right)\to\min \quad , (2.7)$$

$$p_i^k\le\hat{p}_j^k+c_{ij}^k,\ p_i^k\ge 0,\hat{p}_j^k\ge 0\ \left(i\in I,j\in J,k\in K\right) \quad . (2.8)$$

**Теорема 2.1.** *Для того чтобы набор неотрицательных векторов*

$$\left\{X_i^k,Y_j^k,z_{ij}^k,p_i^k,\hat{p}_j^k\,\middle|\,i\in I,j\in J,k\in K\right\}$$

*являлся конкурентным равновесием в модели железнодорожных грузоперевозок, необходимо и достаточно, чтобы набор* $\left\{X_i^k,Y_j^k,z_{ij}^k\,\middle|\,i\in I,j\in J,k\in K\right\}$ *являлся решением экстремальной задачи (2.1)- (2.4), а* $\left\{p_i^k,\hat{p}_j^k\,\middle|\,i\in I,j\in J,k\in K\right\}$ *являлся решением двойственной задачи (2.7)-(2.8). При этом оптимальные значения функционалов в задачах (2.1)-(2.4) и (2.7)-(2.8) равны.*

**Доказательство.** Обозначим множители Лагранжа к ограничениям (2.2) через $p_i^k\ge 0\left(i\in I,k\in K\right)$, а к ограничениям (2.3) - через $\hat{p}_j^k\ge 0\left(j\in J,k\in K\right)$.

Составим функцию Лагранжа задачи (2.1)-(2.4):

$$L\left(X_i^k, Y_j^k, z_{ij}^k, p_i^k, \hat{p}_j^k \middle| i \in I, j \in J, k \in K\right) = \sum_{i \in I} F_i\left(X_i^1, ..., X_i^m\right) -$$

$$-\sum_{j \in J} G_j\left(Y_j^1, ..., Y_j^m\right) - \sum_{i \in I, j \in J, k \in K} c_{ij}^k z_{ij}^k +$$

$$+\sum_{i \in I, k \in K} p_i^k \left(\sum_{j \in J} z_{ij}^k - X_i^k\right) + \sum_{j \in J, k \in K} \hat{p}_j^k \left(Y_j^k - \sum_{i \in I} z_{ij}^k\right) =$$

$$= \sum_{i \in I}\left(F_i\left(X_i^1, ..., X_i^m\right) - \sum_{k \in K} p_i^k X_i^k\right) +$$

$$+\sum_{j \in J}\left(\sum_{k \in K} \hat{p}_j^k Y_j^k - G_j\left(Y_j^1, ..., Y_j^m\right)\right) + \sum_{i \in I, j \in J, k \in K} z_{ij}^k\left(p_i^k - \hat{p}_j^k - c_{ij}^k\right).$$

(2.9)

По теореме Куна-Таккера набор

$$\left\{X_i^k, Y_j^k, z_{ij}^k \middle| i \in I, j \in J, k \in K\right\}$$

является решением задачи (2.1)-(2.4) тогда и только тогда, когда существуют множители Лагранжа $p_i^k \geq 0 \left(i \in I, k \in K\right)$, $\hat{p}_j^k \geq 0 \left(j \in J, k \in K\right)$ такие, что

$$\left\{X_i^k, Y_j^k, z_{ij}^k, p_i^k, \hat{p}_j^k \middle| i \in I, j \in J, k \in K\right\}$$

является конкурентным равновесием в модели железнодорожных грузоперевозок. Из теоремы Куна-Таккера, примененной к задаче (2.7)-(2.8), следует, что множители Лагранжа $p_i^k \geq 0 \left(i \in I, k \in K\right)$, $\hat{p}_j^k \geq 0 \left(j \in J, k \in K\right)$ являются ее решением. Действительно, обозначим множители Лагранжа к ограничениям (2.8) через $z_{ij}^k \geq 0 \left(i \in I, j \in J, k \in K\right)$ и составим функцию Лагранжа задачи (2.7)-(2.8):

$$L\left(p_i^k, \hat{p}_j^k, z_{ij}^k \middle| i \in I, j \in J, k \in K\right) = \sum_{i\in I} \Pi_i\left(p_i^1,\dots,p_i^m\right) +$$

$$+\sum_{j\in J} \pi_j\left(\hat{p}_j^1,\dots,\hat{p}_j^m\right) + \sum_{i\in I, j\in J, k\in K} z_{ij}^k\left(p_i^k - \hat{p}_j^k - c_{ij}^k\right) =$$

$$= \sum_{i\in I}\left(\Pi_i\left(p_i^1,\dots,p_i^m\right) + \sum_{j\in J, k\in K} p_i^k z_{ij}^k\right) +$$

$$+\sum_{j\in J}\left(\pi_j\left(\hat{p}_j^1,\dots,\hat{p}_j^m\right) - \sum_{i\in I, k\in K} \hat{p}_j^k z_{ij}^k\right) - \sum_{i\in I, j\in J, k\in K} z_{ij}^k c_{ij}^k.$$

Полагая $\breve{X}_i^k = \sum_{j\in J} z_{ij}^k \; (i \in I, k \in K)$, $\breve{Y}_j^k = \sum_{i\in I} z_{ij}^k \; (j \in J, k \in K)$, получаем с учетом (2.5), (2.6), что множители Лагранжа $p_i^k \geq 0 \; (i \in I, k \in K)$, $\hat{p}_j^k \geq 0 \; (j \in J, k \in K)$ являются решением задачи (2.7)-(2.8) тогда и только тогда, когда набор $\left\{\breve{X}_i^k, \breve{Y}_j^k, z_{ij}^k, p_i^k, \hat{p}_j^k \middle| i \in I, j \in J, k \in K\right\}$ является конкурентным равновесием в модели железнодорожных грузоперевозок. Равенство функционалов в задачах (2.1)-(2.4) и (2.7)-(2.8) следует из теоремы Фенхеля [12, с. 46-47]. Теорема 2.1 доказана.

Из теоремы 2.1 следует, что конкурентное равновесие в модели железнодорожных перевозок является экономически эффективным.

## 3. Несовершенная конкуренция

Появление на рынке услуг по перевозке грузов частных компаний, инвестировавших средства в обновление вагонного парка, потенциально может привести к повышению тарифов и сокращению объемов перевозок. Для оценки возникающей угрозы рассмотрим задачу, в которой «перевозчик», пользуясь услугами ОАО «РЖД» по перевозке $k$-го товара из $j$-го пункта производства в $i$-й пункт потребления по тарифу ОАО «РЖД» $\tilde{c}_{ij}^k$ и предоставляя принадлежащие ему вагоны, назначает свой тариф на услугу $c_{ij}^k$ в целях максимизации своего дохода

$$\sum_{k\in K, i\in I, j\in J} \left(c_{ij}^k - \tilde{c}_{ij}^k\right) z_{ij}^k(c) \to \max_{c_{ij}^k \geq \tilde{c}_{ij}^k} . \tag{3.1}$$

Здесь функции спроса $z_{ij}^k(\bullet)$ на услуги по перевозке определяются из решения семейства задач (2.1)-(2.4) при различных значениях тарифов $c = \left\{c_{ij}^k \middle| i \in I, j \in J, k \in K\right\}$.

Положим $\Phi(c) = \sum_{k\in K, i\in I, j\in J} \left(c_{ij}^k - \tilde{c}_{ij}^k\right) z_{ij}^k(c)$ и

$$\Lambda\left(c_{ij}^k \middle| i \in I, j \in J, k \in K\right) =$$

$$\min\left\{ \sum_{i\in I} \Pi_i\left(p_i^1, ..., p_i^m\right) + \sum_{j\in J} \pi_j\left(\hat{p}_j^1, ..., \hat{p}_j^m\right) + \right.$$

$$\left. + \sum_{i\in I, j\in J, k\in K} \tau_{ij}^k \left(c_{ij}^k + \hat{p}_j^k - p_i^k\right) \right|$$

$$\left. p_i^k \geq 0, \hat{p}_j^k \geq 0, \tau_{ij}^k \leq 0 \;\; \left(i \in I, j \in J, k \in K\right) \right\}.$$

Функция $\Lambda(c)$ является выпуклой, ее значение равно оптимальному значению функционала задачи (2.1)-(2.4) при значениях параметров $c = \left\{c_{ij}^k \middle| i \in I, j \in J, k \in K\right\}$, равных значениям аргумента функции $\Lambda(c)$. Кроме того,

$$\left(z_{ij}^k(c) \middle| i \in I, j \in J, k \in K\right) \in -\partial\Lambda(c),$$

где $\partial\Lambda(\bullet)$ - субдифференциал функции $\Lambda(\bullet)$.

Таким образом, анализ монопольного извлечения сверхприбыли перевозчиком сводится к анализу двухуровневой задачи оптимизации. В силу выпуклости функции $\Lambda(c)$ имеем

$$\Phi(c)=\sum_{k\in K,i\in I,j\in J}\left(c_{ij}^k-\tilde{c}_{ij}^k\right)z_{ij}^k(c)\leq\Lambda(\tilde{c})-\Lambda(c).$$

Величина $\Phi(c)$ равна сверхприбыли (доходу от инвестиций), которую получает посредник на рынке грузоперевозок, назначающий повышенные тарифы $c=\left\{c_{ij}^k\,\middle|\,i\in I,\,j\in J,k\in K\right\}$ вместо тарифов $\tilde{c}=\left\{\tilde{c}_{ij}^k\,\middle|\,i\in I,\,j\in J,k\in K\right\}$.

. Величина $\Lambda(\tilde{c})-\Lambda(c)$ равна совокупным убыткам потребителей и производителей в результате изменения тарифов. Отношение

$$\theta=\frac{\Phi(c)}{\Lambda(\tilde{c})-\Lambda(c)}$$

является показателем общесистемных потерь в результате изменения тарифов. Если $\theta=1$, то сверхприбыль посредника равна совокупным убыткам потребителей и производителей, и общесистемные потери отсутствуют. Если $\theta<1$, то повышение тарифов приводит не только к перераспределению доходов, но снижает эффективность функционирования экономической системы в целом. Чем меньше $\theta$, тем больше доля общесистемных потерь по сравнению с величиной перераспределяемой прибыли.

### *3.1. Монопольный перевозчик в условиях совершенной конкуренции между производителем и потребителем товаров*

Рассмотрим пример рынка однородного товара, на котором производитель, имеющий функцию себестоимости

$$G(Y)=\begin{cases} BY+AY^2, & \text{если } Y\le \hat{Y},\\ +\infty, & \text{если } Y>\hat{Y},\end{cases}$$

взаимодействует с потребителем, имеющим функцию дохода от используемого груза

$$F(X)=\begin{cases} bX+aX^2, & \text{если } X\le \hat{X},\\ b\hat{X}+a\hat{X}^2, & \text{если } X>\hat{X},\end{cases}$$

и перевозчиком, имеющим себестоимость перевозки единицы груза $\tilde{c}$. Здесь

$$a<0, A>0,\ b>0, B>0, c>0, \hat{X}>0, \hat{Y}>0,\ \hat{X}\le -\frac{b}{2a}.$$

Тогда

$$\Pi(p)=\begin{cases} -\dfrac{(b-p)_+^2}{4a}, & \text{если } b+2a\hat{X}\le p,\\ (b-p)\hat{X}+a\hat{X}^2, & \text{если } b+2a\hat{X}>p,\end{cases}$$

$$\pi(\hat{p})=\begin{cases} \dfrac{(\hat{p}-B)_+^2}{4A}, & \text{если } \hat{p}\le B+2A\hat{Y},\\ (\hat{p}-B)\hat{Y}-A\hat{Y}^2, & \text{если } \hat{p}>B+2A\hat{Y}.\end{cases}$$

Для рассматриваемого примера задача (2.6)-(2.7) имеет вид

$$\Pi(p)+\pi(\hat{p})\to\min\quad, \qquad (3.2)$$

$$p\le \hat{p}+c, p\ge 0, \hat{p}\ge 0. \qquad (3.3)$$

Поскольку функция $\Pi(p)$ не возрастает по переменной $p$, можно положить на оптимальном решении $p=\hat{p}+c$. Тогда задача минимизации (3.2)-(3.3) сводится к следующей задаче:

$$\min\left\{\Pi(\hat{p}+c)+\pi(\hat{p}) \mid \hat{p}\ge 0\right\}.$$

Если потребитель не насыщен и производственные возможности не исчерпаны, т.е.

$$\left(b-B-c\right)_+ < 2\left(A-a\right)\min\left(\hat{X},\hat{Y}\right), \text{ то}$$

$$\frac{\partial\Pi\left(\hat{p}+c\right)}{\partial\hat{p}}+\frac{\partial\pi\left(\hat{p}\right)}{\partial\hat{p}}=0,$$

что эквивалентно

$$\frac{b-\hat{p}-c}{2a}+\frac{\hat{p}-B}{2A}=0.$$

Откуда следует, что

$$\hat{p}=\frac{bA-cA-Ba}{A-a}.$$

Такой цене в пункте производства соответствует объем грузоперевозок

$$z\left(c\right)=\frac{\partial\pi\left(\hat{p}\right)}{\partial\hat{p}}=\frac{\left(\hat{p}-B\right)_+}{2A}=\frac{\left(b-c-B\right)_+}{2\left(A-a\right)}. \qquad (3.4)$$

В рассматриваемом примере задача (3.1) имеет вид

$$\left(c-\tilde{c}\right)z\left(c\right)\to\max_{c\geq\tilde{c}}.$$

Если $0<b-B-\tilde{c}<4\left(A-a\right)\min\left(\hat{X},\hat{Y}\right)$, то решение этой задачи $\breve{c}$ и соответствующий ему объем грузоперевозок $z\left(\breve{c}\right)$ равны

$$\breve{c}=\frac{b-B+\tilde{c}}{2},\ z\left(\breve{c}\right)=\frac{b-B-\tilde{c}}{4\left(A-a\right)}.$$

Подчеркнем, что $z\left(\breve{c}\right)=\dfrac{z\left(\tilde{c}\right)}{2}$, т.е. в рассматриваемом примере появление монопольного посредника (перевозчика) между ОАО «РЖД» и производителем-потребителем товаров приводит к такому повышению тарифов на грузоперевозки, при котором грузооборот уменьшается вдвое. Нетрудно вычислить, что $\theta = \frac{2}{3}$, т.е. извлечение монопольной сверхприбыли посредником сопровождается большим снижением общесистемной эффективности.

### 3.2. Олигополия Курно перевозчиков в условиях совершенной конкуренции между производителем и потребителем товаров

Будем считать выполненными предположения

$$a < 0, A > 0,\ b > 0, B > 0, c > 0, \hat{X} > 0, \hat{Y} > 0,\ \hat{X} \le -\frac{b}{2a},$$

$$0 < b - B - \tilde{c} < 4(A-a)\min\left(\hat{X}, \hat{Y}\right).$$

Из выражения (3.4) можно найти обратные функции спроса на услуги по грузоперевозке, т.е. зависимость между тарифом на перевозки $c$ и спросом на услуги по грузоперевозке $z$ при этом тарифе, в условиях совершенной конкуренции между производителем и потребителем:

$$c(z) = \begin{cases} b - B - 2(A-a)z, \text{ если } 0 \le z \le \dfrac{b-B-\tilde{c}}{2(A-a)}, \\ \tilde{c}, \qquad\qquad \text{если } z > \dfrac{b-B-\tilde{c}}{2(A-a)}. \end{cases} \tag{3.5}$$

Рассмотрим олигополию Курно, в которой выделено $n$ равноправных перевозчиков. Тогда индекс рыночной концентрации Герфиндаля – Хиршмана равен $HHI = \frac{10000}{n}$. Будем считать, что каждый из них, например перевозчик $\xi$, выбирает объем оказываемых услуг по грузоперевозкам $z_\xi$ так, чтобы максимизировать свою прибыль:

$$z_\xi \left( c\left( \sum_{\zeta=1}^{n} z_\zeta \right) - \tilde{c} \right) \to \max_{z_\xi}. \qquad (3.6)$$

Таким образом, прибыль перевозчика $\xi$ зависит не только от его стратегии по выбору объема оказываемых услуг по грузоперевозкам $z_\xi$, но и от стратегий других перевозчиков. В модели Курно предполагается, что компромиссом между интересами перевозчиков является равновесие по Нэшу в игре в нормальной форме, в которой стратегиями перевозчиков являются объёмы грузоперевозок $z_\xi \geq 0$, а функциями выигрыша - их прибыли (3.6). В случае, когда обратные функции спроса имеют вид (3.5), в этой игре существует единственное равновесие по Нэшу (см., например, [13,14]), и соответствующий суммарный объём грузоперевозок $z$ равен

$$z = \frac{n}{2(n+1)} \frac{(b - B - \tilde{c})}{(A - a)}. \qquad (3.7)$$

Определим показатель, характеризующий общесистемные потери при олигополистическом росте тарифов, как

$$\theta = \frac{(c(z) - \tilde{c})z}{\Lambda(\tilde{c}) - \Lambda(c(z))}.$$

Нетрудно вычислить, что в случае олигополии Курно $\theta = \dfrac{2n}{2n+1}$. Из (3.7) видно, что увеличение числа перевозчиков приводит к увеличению объёма грузоперевозок и снижению тарифов на грузоперевозки (к монотонности обратных функций спроса). В случаях, когда $n = 1$, объём грузоперевозок, определяемый по формуле (3.7), равен объёму грузоперевозок $z(\breve{c})$ монопольного грузоперевозчика. При $n \to +\infty$ объём грузоперевозок стремится к объему перевозок в условиях совершенной конкуренции: показатель $\theta \to 1$, т.е. общесистемные потери исчезают. В зависимости от индекса рыночной концентрации Герфиндаля – Хиршмана $HHI$ грузооборот уменьшается в $\frac{10000}{10000+HHI}$ раз, а показатель общесистемных потерь $\theta$ в $\frac{20000}{20000+HHI}$ раз.

### 3.3. Иерархия равновесий по Штаккельбергу

Будем считать выполненными предположения

$$a<0, A>0,\ b>0, B>0, c>0, \hat{X}>0, \hat{Y}>0,\ \hat{X} \le -\frac{b}{2a},$$

$$0<b-B-\tilde{c}<4(A-a)\min\left(\hat{X},\hat{Y}\right).$$

$$0<b-B-\tilde{c}<4(A-a)\min\left(\hat{X},\hat{Y}\right) \text{ и } b+B-\tilde{c} \le 4\hat{X}\left(A-2a\right).$$

Рассмотрим пример рынка однородного товара, на котором потребитель, управляя заказом на перевозку продукции, максимизирует свою прибыль от сбыта продукции, принимая во внимание цены на продукцию у производителя, тариф перевозчика, прибыль от использования товара $F(X)$, описанную в разд. 3.1.

В результате решения задачи потребителя определяется функция, описывающая зависимость объема заказа на грузоперевозку между потребителем и производителем от цен на продукцию у производителя и тарифа перевозчика:

$$z\left(\hat{p}+c\right)=\begin{cases} -\dfrac{\left(b-\hat{p}-c\right)}{2a}, & \text{если } b+2a\hat{X} \le \hat{p}+c, \\ \hat{X}, & \text{если } b+2a\hat{X} > \hat{p}+c. \end{cases} \tag{3.8}$$

Будем считать, что производитель, управляя ценами на производимый им товар, максимизирует свою прибыль, принимая во внимание то, как формируется спрос на его продукцию, т.е. результирующую функцию спроса потребителя (3.8) и себестоимость производства $G(Y)$, описанную в разд. 3.1. В результате он назначает цену $\hat{p}$, которая является решением задачи

$$\hat{p}z(\hat{p}+c)-G(z(\hat{p}+c))\to\max_{\hat{p}\geq 0} \tag{3.9}$$

Из решения задачи (3.9) получаем, что если $b+B-2\hat{X}(A-2a)\leq c$, то

$$\hat{p}=\frac{(b-c)(A-a)+aB}{A-2a}.$$

Подставляя найденное выражение в (3.8), получаем функцию спроса на грузоперевозки в зависимости от тарифа $c$:

$$h(c)=\begin{cases}\dfrac{b+B-c}{2(A-2a)}, & \text{если } b+B-2\hat{X}(A-2a)\leq c,\\ \hat{X}, & \text{если } b+B-2\hat{X}(A-2a)>c.\end{cases} \tag{3.10}$$

Перевозчик, управляя тарифом $c$, максимизирует свою прибыль $(c-\tilde{c})h(c)$, принимая во внимание себестоимость услуги по грузоперевозке $\tilde{c}$ (тариф ОАО «РЖД») инфраструктуры и функцию зависимости спроса $h(c)$ на перевозки от тарифа на грузоперевозку $c$. Решая эту задачу с функцией (3.10), получаем, что значение тарифа

$$\widehat{c}=\frac{b+B+\tilde{c}}{2}$$

больше тарифа $\breve{c}$ в случае монопольного перевозчика в условиях совершенной конкуренции между производителем и потребителем товаров из разд. 3.1. Тарифу $\hat{c}$ соответствует объём грузоперевозок

$$h(\hat{c}) = \frac{b + B - \tilde{c}}{4(A - 2a)}.$$

Подчеркнём, что во всех рассмотренных случаях увеличение тарифов за услуги ОАО «РЖД» уменьшает объёмы грузоперевозок.

# 4. Конкурентное равновесие в модели железнодорожных грузоперевозок с коммуникационными ограничениями. Два подхода к формированию тарифов

Будем считать, что перевозка $k$-го товара из $j$-го пункта производства в $i$-й пункт потребления осуществляется по маршруту $\Gamma_{ij}^k$, который описывается множеством ребер графа, соединяющего соседние узловые станции, через которые проходит маршрут. Обозначим через $\tilde{t}_{\alpha\beta}$ себестоимость для ОАО «РЖД» перевозки единицы груза между соседними узловыми станциями со станции $\alpha$ на станцию $\beta$ и через $\tilde{\lambda}_\alpha$ - себестоимость для ОАО «РЖД» приема и отправления единицы груза через узловую станцию $\alpha$. Тогда себестоимость для ОАО «РЖД» перевозки единицы $k$-го товара из $j$-го пункта производства в $i$-й пункт потребления составит

$$\tilde{c}_{ij}^k = \sum_{(\alpha,\beta)\in\Gamma_{ij}^k} \tilde{t}_{\alpha\beta} + \sum_{(\alpha,\beta)\in\Gamma_{ij}^k} \tilde{\lambda}_\alpha .$$

Обозначим пропускную способность железной дороги по перевозке грузов из узловой станции $\alpha$ на узловую станцию $\beta$

через $V_{\alpha\beta}$. Пропускную способность узловой станции $\alpha$ по приёму и отправлению грузов обозначим $M_\alpha$.

Рассмотрим задачу о максимизации прибыли экономической системы с учётом коммуникационных ограничений:

$$\sum_{i\in I} F_i\left(X_i^1,\dots,X_i^m\right)-\sum_{j\in J} G_j\left(Y_j^1,\dots,Y_j^m\right)-$$

$$-\sum_{i\in I,\,j\in J,\,k\in K}\left(\sum_{(\alpha,\beta)\in\Gamma_{ij}^k}\tilde{t}_{\alpha\beta}+\sum_{(\alpha,\beta)\in\Gamma_{ij}^k}\tilde{\lambda}_\alpha\right)z_{ij}^k\to\max \qquad , \qquad (4.1)$$

$$X_i^k\le\sum_{j\in J} z_{ij}^k \quad \left(i\in I, k\in K\right), \qquad (4.2)$$

$$Y_j^k\ge\sum_{i\in I} z_{ij}^k \quad \left(j\in J, k\in K\right), \qquad (4.3)$$

$$\sum_{k\in K}\ \sum_{\left\{i,j\,\middle|\,(\alpha,\beta)\in\Gamma_{ij}^k\right\}} z_{ij}^k\le V_{\alpha\beta} \qquad \textit{для любых } \left(\alpha,\beta\right) \qquad , \qquad (4.4)$$

$$\sum_{k\in K}\ \sum_{\left\{i,j\,\middle|\,(\alpha,\beta)\in\Gamma_{ij}^k\right\}} z_{ij}^k\le M_\alpha \qquad \textit{для любого } \alpha \qquad , \qquad (4.5)$$

$$z_{ij}^k\ge 0 \qquad \left(i\in I, j\in J, k\in K\right). \qquad (4.6)$$

Отметим, что при решении задачи (4.1)-(4.6) экстремум следует искать также за счет выбора маршрутов перевозок $\Gamma_{ij}^k$.

Двойственной по Фенхелю [12, с. 46-47] экстремальной задачей к задаче (4.1)-(4.6) является задача выпуклого программирования:

$$\sum_{i\in I}\Pi_i\left(p_i^1,\dots,p_i^m\right)+\sum_{j\in J}\pi_j\left(\hat{p}_j^1,\dots,\hat{p}_j^m\right)+$$

$$+\sum_{(\alpha,\beta)\in\Gamma_{ij}^k} t_{\alpha\beta}V_{\alpha\beta}+\sum_{(\alpha,\beta)\in\Gamma_{ij}^k}\lambda_\alpha M_\alpha\to\min \qquad , (4.7)$$

$$p_i^k \le \hat{p}_j^k + \sum_{(\alpha,\beta)\in\Gamma_{ij}^k} \left(t_{\alpha\beta} + \tilde{t}_{\alpha\beta}\right) + \sum_{(\alpha,\beta)\in\Gamma_{ij}^k} \left(\lambda_\alpha + \tilde{\lambda}_\alpha\right), \qquad (4.8)$$

$$p_i^k \ge 0, \hat{p}_j^k \ge 0, t_{\alpha\beta} \ge 0,$$

$$\lambda_\alpha \ge 0 \quad \left(i \in I, j \in J, k \in K, (\alpha,\beta) \in \Gamma_{ij}^k\right) \qquad . (4.9)$$

Здесь множители Лагранжа $t_{\alpha\beta}$ к коммуникационным ограничениям (4.4) интерпретируются как теневые наценки при эксплуатации ограниченной пропускной способности железной дороги, соединяющей узловые станции $\alpha$ и $\beta$, а множители Лагранжа к коммуникационным ограничениям (4.5) $\lambda_\alpha$ интерпретируются как теневые наценки при эксплуатации ограниченной пропускной способности узловой станции $\alpha$.

**Определение.** *Будем говорить, что набор неотрицательных чисел*

$$\left\{X_i^k, Y_j^k, z_{ij}^k, p_i^k, \hat{p}_j^k, t_{\alpha\beta}, \lambda_\alpha \,\middle|\, i \in I, j \in J, k \in K, (\alpha,\beta) \in \Gamma_{ij}^k\right\}$$

*является конкурентным равновесием в модели железнодорожных грузоперевозок с коммуникационными ограничениями, если*

*1) для любого* $i \in I$

$$\left(X_i^1,...,X_i^m\right) \in Arg\max\left\{F_i\left(x_i^1,...,x_i^m\right) - \sum_{i=1}^{k} p_i^k x_i^k \,\middle|\, \left(x_i^1,...,x_i^m\right) \ge 0\right\},$$

*2) для любого* $j \in J$

$$\left(Y_j^1,...,Y_j^m\right) \in Arg\max\left\{\sum_{k=1}^{m} \hat{p}_j^k y_j^k - G_j\left(y_j^1,...,y_j^m\right) \,\middle|\, \left(y_j^1,...,y_j^m\right) \ge 0\right\},$$

*3) для любых* $i \in I$, $k \in K$ $X_i^k \le \sum_{j\in J} z_{ij}^k, \; p_i^k\left(X_i^k - \sum_{j\in J} z_{ij}^k\right) = 0,$

4) *для любых* $j \in J$ , $k \in K$ $Y_j^k \geq \sum_{i \in I} z_{ij}^k$, $\hat{p}_j^k \left( Y_j^k - \sum_{i \in I} z_{ij}^k \right) = 0$,

5) *для любых* $i \in I$ , $j \in J$ , $k \in K$

$$p_i^k \leq \hat{p}_j^k + \sum_{(\alpha,\beta) \in \Gamma_{ij}^k} \left( t_{\alpha\beta} + \tilde{t}_{\alpha\beta} \right) + \sum_{(\alpha,\beta) \in \Gamma_{ij}^k} \left( \lambda_\alpha + \tilde{\lambda}_\alpha \right),$$

$$z_{ij}^k \left( \hat{p}_j^k + \sum_{(\alpha,\beta) \in \Gamma_{ij}^k} \left( t_{\alpha\beta} + \tilde{t}_{\alpha\beta} \right) + \sum_{(\alpha,\beta) \in \Gamma_{ij}^k} \left( \lambda_\alpha + \tilde{\lambda}_\alpha \right) - p_i^k \right) = 0,$$

$$z_{ij}^k \geq 0,$$

6) *для любых* $(\alpha, \beta)$

$$\sum_{k \in K} \sum_{\left\{ i,j \middle| (\alpha,\beta) \in \Gamma_{ij}^k \right\}} z_{ij}^k \leq V_{\alpha\beta},\ t_{\alpha\beta} \left( \sum_{k \in K} \sum_{\left\{ i,j \middle| (\alpha,\beta) \in \Gamma_{ij}^k \right\}} z_{ij}^k - V_{\alpha\beta} \right) = 0,$$

7) *для любых* $(\alpha, \beta)$

$$\sum_{k \in K} \sum_{\left\{ i,j \middle| (\alpha,\beta) \in \Gamma_{ij}^k \right\}} z_{ij}^k \leq M_\alpha,\ \lambda_\alpha \left( \sum_{k \in K} \sum_{\left\{ i,j \middle| (\alpha,\beta) \in \Gamma_{ij}^k \right\}} z_{ij}^k - M_\alpha \right) = 0.$$

**Теорема 4.1.** *Для того чтобы набор неотрицательных векторов*

$$\left\{ X_i^k, Y_j^k, z_{ij}^k, p_i^k, \hat{p}_j^k, t_{\alpha\beta}, \lambda_\alpha \middle| i \in I, j \in J, k \in K, (\alpha, \beta) \in \Gamma_{ij}^k \right\} \quad (4.10)$$

*являлся конкурентным равновесием в модели железнодорожных грузоперевозок с коммуникационными ограничениями, необходимо и достаточно, чтобы набор*

$\left\{ X_i^k, Y_j^k, z_{ij}^k \middle| i \in I, j \in J, k \in K \right\}$ *являлся решением экстремальной задачи (4.1)- (4.6), а*

$$\left\{ p_i^k, \hat{p}_j^k, t_{\alpha\beta}, \lambda_\alpha \middle| i \in I, j \in J, k \in K, (\alpha,\beta) \in \Gamma_{ij}^k \right\}$$

*являлся решением двойственной задачи (4.7)-(4.9)*[3]*. При этом оптимальные значения функционалов в задачах (4.1)-(4.6) и (4.7)-(4.9) равны.*

**Доказательство.** Обозначим множители Лагранжа к ограничениям (4.2) через $p_i^k \geq 0 \ (i \in I, k \in K)$, а к ограничениям (4.3) - через $\hat{p}_j^k \geq 0 \ (j \in J, k \in K)$, к ограничениям (4.4) - через $t_{\alpha\beta} \geq 0 \ \left((\alpha,\beta) \in \Gamma_{ij}^k\right)$, к ограничениям (4.5) - через $\lambda_\alpha \geq 0 \ \left((\alpha,\beta) \in \Gamma_{ij}^k\right)$. Составим функцию Лагранжа задачи (2.1)-(2.4):

$$L\left( X_i^k, Y_j^k, z_{ij}^k, p_i^k, \hat{p}_j^k, t_{\alpha\beta}, \lambda_\alpha \middle| i \in I, j \in J, k \in K, (\alpha,\beta) \in \Gamma_{ij}^k \right) =$$

$$= \sum_{i \in I} F_i\left(X_i^1, \ldots, X_i^m\right) - \sum_{j \in J} G_j\left(Y_j^1 \ldots, Y_j^m\right) -$$

$$- \sum_{i \in I, j \in J, k \in K} \left( \sum_{(\alpha,\beta) \in \Gamma_{ij}^k} \tilde{t}_{\alpha\beta} + \sum_{(\alpha,\beta) \in \Gamma_{ij}^k} \tilde{\lambda}_\alpha \right) z_{ij}^k +$$

$$+ \sum_{i \in I, k \in K} p_i^k \left( \sum_{j \in J} z_{ij}^k - X_i^k \right) + \sum_{j \in J, k \in K} \tilde{p}_j^k \left( Y_j^k - \sum_{i \in I} z_{ij}^k \right) +$$

---

[3] Будем говорить в этом случае, что конкурентное равновесие (4.10) соответствует паре взаимно двойственных задач выпуклого программирования (4.1)-(4.6) и (4.7)-(4.9).

$$+\sum_{(\alpha,\beta)} t_{\alpha\beta}\left(V_{\alpha\beta}-\sum_{\left\{(i.j),k\,\middle|\,(\alpha,\beta)\in\Gamma_{ij}^{k}\right\}} z_{ij}^{k}\right)+$$

$$+\sum_{\alpha}\lambda_{\alpha}\left(M_{\alpha}-\sum_{\left\{(i.j),k\,\middle|\,(\alpha,\beta)\in\Gamma_{ij}^{k}\right\}} z_{ij}^{k}\right)=$$

$$=\sum_{i\in I}\left(F_i\left(X_i^1,\ldots,X_i^m\right)-\sum_{k\in K}p_i^k X_i^k\right)+$$

$$+\sum_{j\in J}\left(\sum_{k\in K}\hat{p}_j^k Y_j^k-G_j\left(Y_j^1,\ldots,Y_j^m\right)\right)+$$

$$+\sum_{i\in I,\,j\in J,k\in K} z_{ij}^{k}\left(p_i^k-\hat{p}_j^k-\sum_{(\alpha,\beta)\in\Gamma_{ij}^{k}}\left(t_{\alpha\beta}+\tilde{t}_{\alpha\beta}\right)-\sum_{(\alpha,\beta)\in\Gamma_{ij}^{k}}\left(\lambda_{\alpha}+\tilde{\lambda}_{\alpha}\right)\right).$$

Если маршруты грузоперевозок между пунктами производства и потребления $\left\{\Gamma_{ij}^{k}\,\middle|\,i\in I,\,j\in J,k\in K\right\}$ зафиксированы, то по теореме Куна-Таккера набор $\left\{X_i^k,Y_j^k,z_{ij}^k\,\middle|\,i\in I,\,j\in J,k\in K\right\}$ является решением задачи (4.1)-(4.6) тогда и только тогда, когда существуют множители Лагранжа

$$p_i^k\geq 0\ \left(i\in I,k\in K\right),\ \hat{p}_j^k\geq 0\ \left(j\in J,k\in K\right),$$

$$t_{\alpha\beta}\geq 0\quad\left((\alpha,\beta)\in\Gamma_{ij}^{k}\right),\ \lambda_{\alpha}\geq 0\quad\left((\alpha,\beta)\in\Gamma_{ij}^{k}\right)$$

такие, что

$$\left\{X_i^k,Y_j^k,z_{ij}^k,p_i^k,\hat{p}_j^k,t_{\alpha\beta},\lambda_{\alpha}\,\middle|\,i\in I,\,j\in J,k\in K,(\alpha,\beta)\in\Gamma_{ij}^{k}\right\}$$

удовлетворяет условиям 1-7 из определения конкурентного равновесия в модели железнодорожных грузоперевозок с коммуникационными ограничениями.

Из теоремы Куна-Таккера, примененной к задаче (4.7)-(4.9), следует, что множители Лагранжа

$$p_i^k \geq 0 \; \left(i \in I, k \in K\right), \; \hat{p}_j^k \geq 0 \; \left(j \in J, k \in K\right),$$

$$t_{\alpha\beta} \geq 0 \quad \left(\left(\alpha, \beta\right) \in \Gamma_{ij}^k\right), \; \lambda_\alpha \geq 0 \quad \left(\left(\alpha, \beta\right) \in \Gamma_{ij}^k\right)$$

являются её решением. Действительно, обозначим множители Лагранжа к ограничениям (4.8) через $z_{ij}^k \geq 0 \; \left(i \in I, j \in J, k \in K\right)$ и составим функцию Лагранжа задачи (4.7)-(4.9):

$$L\left(p_i^k, \hat{p}_j^k, z_{ij}^k, t_{\alpha\beta} \geq 0, \lambda_\alpha \middle| i \in I, j \in J, k \in K, (\alpha, \beta) \in \Gamma_{ij}^k\right) =$$

$$= \sum_{i \in I} \Pi_i \left(p_i^1, \ldots, p_i^m\right) + \sum_{j \in J} \pi_j \left(\hat{p}_j^1, \ldots, \hat{p}_j^m\right) + \sum_{(\alpha, \beta) \in \Gamma_{ij}^k} t_{\alpha\beta} V_{\alpha\beta} + \sum_{(\alpha, \beta) \in \Gamma_{ij}^k} \lambda_\alpha M_\alpha +$$

$$+ \sum_{i \in I, j \in J, k \in K} z_{ij}^k \left( p_i^k - \hat{p}_j^k - \sum_{(\alpha, \beta) \in \Gamma_{ij}^k} \left(t_{\alpha\beta} + \tilde{t}_{\alpha\beta}\right) - \sum_{(\alpha, \beta) \in \Gamma_{ij}^k} \left(\lambda_\alpha + \tilde{\lambda}_\alpha\right) \right) =$$

$$= - \sum_{i \in I, j \in J, k \in K} z_{ij}^k \left( \sum_{(\alpha, \beta) \in \Gamma_{ij}^k} \tilde{t}_{\alpha\beta} + \sum_{(\alpha, \beta) \in \Gamma_{ij}^k} \tilde{\lambda}_\alpha \right) +$$

$$+ \sum_{i \in I} \left( \Pi_i \left(p_i^1, \ldots, p_i^m\right) + \sum_{j \in J, k \in K} p_i^k z_{ij}^k \right) +$$

$$+ \sum_{j \in J} \left( \pi_j \left(\hat{p}_j^1, \ldots, \hat{p}_j^m\right) - \sum_{i \in I, k \in K} \hat{p}_j^k z_{ij}^k \right) +$$

$$+\sum_{(\alpha,\beta)\in\Gamma_{ij}^k} t_{\alpha\beta}\left(V_{\alpha\beta}-\sum_{k\in K}\sum_{\left\{i,j\,\middle|\,(\alpha,\beta)\in\Gamma_{ij}^k\right\}} z_{ij}^k\right)+$$

$$+\sum_{(\alpha,\beta)\in\Gamma_{ij}^k} \lambda_{\alpha}\left(M_{\alpha}-\sum_{k\in K}\sum_{\left\{i,j\,\middle|\,(\alpha,\beta)\in\Gamma_{ij}^k\right\}} z_{ij}^k\right).$$

Полагая

$$\breve{X}_i^k=\sum_{j\in J} z_{ij}^k \quad (i\in I, k\in K),\ \breve{Y}_j^k=\sum_{i\in I} z_{ij}^k \quad (j\in J, k\in K),$$

получаем с учётом (2.5), (2.6), что если маршруты грузоперевозок между пунктами производства и потребления $\left\{\Gamma_{ij}^k \middle| i\in I, j\in J, k\in K\right\}$ зафиксированы, то множители Лагранжа

$$p_i^k\ge 0\ (i\in I, k\in K),\ \hat{p}_j^k\ge 0\ (j\in J, k\in K),$$

$$t_{\alpha\beta}\ge 0\quad \left((\alpha,\beta)\in\Gamma_{ij}^k\right),\ \lambda_{\alpha}\ge 0\quad \left((\alpha,\beta)\in\Gamma_{ij}^k\right)$$

являются решением задачи (4.7)-(4.9) тогда и только тогда, когда набор

$$\left\{\breve{X}_i^k, \breve{Y}_j^k, z_{ij}^k, p_i^k, \hat{p}_j^k, t_{\alpha\beta}, \lambda_{\alpha} \middle| i\in I, j\in J, k\in K, (\alpha,\beta)\in\Gamma_{ij}^k\right\}$$

удовлетворяет условиям 1-7 из определения конкурентного равновесия в модели железнодорожных грузоперевозок с коммуникационными ограничениями. Равенство функционалов в задачах (4.1)-(4.6) и (4.7)-(4.9) следует из теоремы Фенхеля [12, с,46-47]. Теорема 4.1 доказана.

Конкурентное равновесие в модели железнодорожных грузоперевозок с коммуникационными ограничениями является эффективным с экономической точки зрения, если считать, что маршруты грузоперевозок $\left\{\Gamma_{ij}^{k} \middle| i \in I, j \in J, k \in K\right\}$ зафиксированы. Однако при решении задачи (4.1)-(4.6) увеличить экономический результат рассматриваемой системы можно также за счет выбора маршрутов грузоперевозок $\left\{\Gamma_{ij}^{k} \middle| i \in I, j \in J, k \in K\right\}$.

**Теорема 4.2.** *Для того чтобы конкурентное равновесие в модели железнодорожных грузоперевозок с коммуникационными ограничениями*

$$\left\{X_i^k, Y_j^k, z_{ij}^k, p_i^k, \hat{p}_j^k, t_{\alpha\beta}, \lambda_\alpha \middle| i \in I, j \in J, k \in K, (\alpha,\beta) \in \Gamma_{ij}^k\right\}$$

*соответствовало паре взаимно двойственных задач (4.1)-(4.6) и (4.7)-(4.9) с оптимизацией в задаче (4.1)-(4.6) по выбору маршрутов грузоперевозок* $\left\{\Gamma_{ij}^{k} \middle| i \in I, j \in J, k \in K\right\}$, *необходимо, чтобы:*

- *для любых пунктов производства* $j$ *и потребления* $i$, *для которых существует грузоперевозка* $z_{ij}^k > 0$ *товара* $k \in K$ *выполнялось условие* 8;
- *для любого маршрута грузоперевозки* $\tilde{\Gamma}_{ij}$ *из пункта производства* $j$ *в пункт потребления* $i$ *выполнялось неравенство*

$$\sum_{(\alpha,\beta)\in\Gamma_{ij}^k} \left(t_{\alpha\beta} + \tilde{t}_{\alpha\beta}\right) + \sum_{(\alpha,\beta)\in\Gamma_{ij}^k} \left(\lambda_\alpha + \tilde{\lambda}_\alpha\right) \le$$

$$\le \sum_{(\tilde{\alpha},\tilde{\beta})\in\tilde{\Gamma}_{ij}} \left(t_{\tilde{\alpha}\tilde{\beta}} + \tilde{t}_{\tilde{\alpha}\tilde{\beta}}\right) + \sum_{(\tilde{\alpha},\tilde{\beta})\in\tilde{\Gamma}_{ij}} \left(\lambda_{\tilde{\alpha}} + \tilde{\lambda}_{\tilde{\alpha}}\right) \quad ,$$

*и достаточно, чтобы для любых маршрутов грузоперевозок* $\left\{\Gamma_{ij}^{k} \middle| i \in I, j \in J, k \in K\right\}$ *выполнялось условие 8.*

**Доказательство.** Зафиксируем значения двойственных переменных $\left\{t_{\alpha\beta},\lambda_\alpha\right\}$ и положим

$$c_{ij}=\min_{\tilde{\Gamma}_{ij}}\sum_{(\tilde{\alpha},\tilde{\beta})\in\tilde{\Gamma}_{ij}}\left(t_{\tilde{\alpha}\tilde{\beta}}+\tilde{t}_{\tilde{\alpha}\tilde{\beta}}\right)+\sum_{(\tilde{\alpha},\tilde{\beta})\in\tilde{\Gamma}_{ij}}\left(\lambda_{\tilde{\alpha}}+\tilde{\lambda}_{\tilde{\alpha}}\right),$$

где $\tilde{\Gamma}_{ij}$ - произвольный возможный маршрут перевозки грузов из пункта производства $j$ в пункт потребления $i$. Сопоставим задаче оптимизации по переменным $\left\{p_i^k,\hat{p}_j^k \middle| i\in I, j\in J, k\in K\right\}$ (4.7)-(4.9) задачу (2.7)-(2.8), значение функционала в которой по теореме 2.1 равно оптимальному значению функционала в задаче (2.1)-(2.4). Тогда утверждение теоремы 4.2 следует из замечания о том, что оптимальное значение функционала в задаче (2.1)-(2.4) монотонно не возрастает по параметрам $\left\{c_{ij}\middle| i\in I, j\in J\right\}$. Обозначим через $z_{ij}^k\left(\tilde{\Gamma}_{ij}\right)$ объём перевозок $k$-го товара из $j$-го пункта производства в $i$-й пункт потребления по маршруту $\tilde{\Gamma}_{ij}$ и рассмотрим вспомогательную задачу:

$$\sum_{i\in I}F_i\left(X_i^1,\dots,X_i^m\right)-\sum_{j\in J}G_j\left(Y_j^1,\dots,Y_j^m\right)-$$

$$-\sum_{i\in I,j\in J,k\in K}\sum_{\tilde{\Gamma}_{ij}}\left(\sum_{(\alpha,\beta)\in\tilde{\Gamma}_{ij}}\tilde{t}_{\alpha\beta}+\sum_{(\alpha,\beta)\in\tilde{\Gamma}_{ij}}\tilde{\lambda}_\alpha\right)z_{ij}^k\left(\tilde{\Gamma}_{ij}\right)\to\max \quad (4.1')$$

$$X_i^k\le\sum_{j\in J}\sum_{\tilde{\Gamma}_{ij}}z_{ij}^k\left(\tilde{\Gamma}_{ij}\right)\quad\left(i\in I,k\in K\right), \quad (4.2')$$

$$Y_j^k\ge\sum_{i\in I}\sum_{\tilde{\Gamma}_{ij}}z_{ij}^k\left(\tilde{\Gamma}_{ij}\right)\quad\left(j\in J,k\in K\right), \quad (4.3')$$

$$\sum_{k\in K}\sum_{\left\{\tilde{\Gamma}_{ij}\middle|(\alpha,\beta)\in\tilde{\Gamma}_{ij}\right\}}z_{ij}^k\left(\tilde{\Gamma}_{ij}\right)\le V_{\alpha\beta}\ \text{ для любых } (\alpha,\beta), \quad (4.4')$$

$$\sum_{k\in K}\sum_{\{\tilde{\Gamma}_{ij}|(\alpha,\beta)\in\tilde{\Gamma}_{ij}\}} z_{ij}^k\left(\tilde{\Gamma}_{ij}\right)\le M_\alpha \quad \textit{для любого } \alpha, \qquad (4.5')$$

$$z_{ij}^k\left(\tilde{\Gamma}_{ij}\right)\ge 0 \qquad \left(i\in I,\, j\in J, k\in K, \tilde{\Gamma}_{ij}\right). \qquad (4.6')$$

Полагая

$$z_{ij}^k\left(\tilde{\Gamma}_{ij}\right)=\begin{cases} z_{ij}^k, \ \textit{если } \tilde{\Gamma}_{ij}=\Gamma_{ij}^k, \\ 0, \ \textit{если } \tilde{\Gamma}_{ij}=\Gamma_{ij}^k, \end{cases}$$

получаем, что

$$\left\{X_i^k, Y_j^k, z_{ij}^k\left(\tilde{\Gamma}_{ij}\right) \middle| i\in I, j\in J, k\in K, \tilde{\Gamma}_{ij}\right\}$$

является допустимым решением задачи (4.1’)-(4.6’). Если множители Лагранжа, являющиеся решением задачи (4.7)-(4.9), удовлетворяют условию 8, то по теореме Куна-Таккера выполняются достаточные условия оптимальности решения. Если же $z_{ij}^k > 0$ и условие 8 нарушается для некоторого альтернативного маршрута $\tilde{\Gamma}_{ij}$, то зафиксируем значения двойственных переменных $\{t_{\alpha\beta}, \lambda_\alpha\}$ и положим

$$c_{ij}=\min_{\tilde{\Gamma}_{ij}}\sum_{(\tilde{\alpha},\tilde{\beta})\in\tilde{\Gamma}_{ij}}\left(t_{\tilde{\alpha}\tilde{\beta}}+\tilde{t}_{\tilde{\alpha}\tilde{\beta}}\right)+\sum_{(\tilde{\alpha},\tilde{\beta})\in\tilde{\Gamma}_{ij}}\left(\lambda_{\tilde{\alpha}}+\tilde{\lambda}_{\tilde{\alpha}}\right).$$

Сопоставим задаче оптимизации по переменным $\left\{p_i^k, \hat{p}_j^k \middle| i\in I, j\in J, k\in K\right\}$ (4.7)-(4.9) задачу (2.7)-(2.8), значение функционала в которой по теореме 2.1 равно оптимальному значению функционала в задаче (2.1)-(2.4). Заметим, что оптимальное значение функционала в задаче (2.1)-(2.4) монотонно не возрастает по параметрам $\left\{c_{ij} \middle| i\in I, j\in J\right\}$. Если же $z_{ij}^k > 0$, то изменение маршрута перевозки $k$-го товара из $j$-го пункта производства в $i$-й пункт потребления с маршрута $\Gamma_{ij}^k$ на маршрут $\tilde{\Gamma}_{ij}$ увеличивает значение функционала в задаче (2.1)-(2.4), а значит, и в задаче (4.1')-(4.6'), т.е. $\left\{X_i^k, Y_j^k, z_{ij}^k\left(\tilde{\Gamma}_{ij}\right) \middle| i \in I, j \in J, k \in K, \tilde{\Gamma}_{ij}\right\}$ не является оптимальным решением. Теорема 4.2. доказана.

Поскольку условие 8 не зависит от типа перевозимого товара, в модели железнодорожных грузоперевозок с коммуникационными ограничениями максимизации прибыли может явиться специализация на перевозке высокодоходных грузов в ущерб низко прибыльным, но социально значимым грузам. Являясь государственной компанией, ОАО «РЖД» несет ответственность за обеспечение перевозки различных, в том числе низко прибыльных грузов. Возможны два подхода к организации грузоперевозок низкой доходности: прямое субсидирование и перекрёстное субсидирование. Как рассчитать тарифы, обеспечивающие за счёт перекрёстного субсидирования перевозку низко прибыльных грузов? Для ответа на этот вопрос рассмотрим сначала вспомогательную экстремальную задачу максимизации прибыли (4.1), к ограничениям которой (4.2)-(4.5) добавим ограничения на перевозку низко прибыльных грузов:

$$z_{ij}^k \geq v_{ij}^k \qquad \left(i \in I, j \in J, k \in K\right). \tag{4.10}$$

Здесь правые части ограничений (4.10) $v_{ij}^k \geq 0$ позволяют учесть требования к перевозке различных грузов.

Двойственная задача к задаче (4.1)-(4.5), (4.10) имеет вид

$$\sum_{i \in I} \Pi_i\left(p_i^1, \ldots, p_i^m\right) + \sum_{j \in J} \pi_j\left(\hat{p}_j^1, \ldots, \hat{p}_j^m\right) + \sum_{(\alpha,\beta) \in \Gamma_{ij}^k} t_{\alpha\beta} V_{\alpha\beta} +$$

$$+ \sum_{(\alpha,\beta) \in \Gamma_{ij}^k} \lambda_\alpha M_\alpha - \sum_{i \in I, j \in J} \gamma_{ij}^k v_{ij}^k \to \min \quad , \tag{4.11}$$

$$p_i^k + \gamma_{ij}^k = \hat{p}_j^k + \sum_{(\alpha,\beta) \in \Gamma_{ij}^k} \left(t_{\alpha\beta} + \tilde{t}_{\alpha\beta}\right) + \sum_{(\alpha,\beta) \in \Gamma_{ij}^k} \left(\lambda_\alpha + \tilde{\lambda}_\alpha\right), \tag{4.12}$$

$$p_i^k \geq 0, \hat{p}_j^k \geq 0, t_{\alpha\beta} \geq 0,$$

$$\lambda_\alpha \ge 0,\ \gamma_{ij}^k \ge 0 \ \left(i \in I, j \in J, k \in K, (\alpha,\beta) \in \Gamma_{ij}^k\right). \qquad (4.13)$$

Множители Лагранжа $\gamma_{ij}^k \ (i \in I, j \in J, k \in K)$ к ограничениям (4.10) интерпретируются как необходимые дотации за перевозку $k$-го товара из $j$-го пункта производства в $i$-й пункт потребления. Дотационный механизм сложен для эффективного администрирования, поэтому предпочтительнее предусмотреть в тарифах перекрёстное субсидирование грузоперевозок. Для расчёта соответствующих тарифов можно предложить механизм фиктивного разыгрывания. Решим сначала вспомогательную задачу (4.1)-(4.5), (4.10) и рассчитаем квоты на использование ограниченных ресурсов по пропускной способности

$$V_{\alpha\beta}^k = \sum_{\{i,j \mid (\alpha,\beta)\in\Gamma_{ij}^k\}} z_{ij}^k,$$

$$M_\alpha^k = \sum_{\{i,j \mid (\alpha,\beta)\in\Gamma_{ij}^k\}} z_{ij}^k.$$

Рассмотрим задачу (4.1)-(4.3), (4.6) с дополнительными ограничениями:

$$V_{\alpha\beta}^k \ge \sum_{\{i,j \mid (\alpha,\beta)\in\Gamma_{ij}^k\}} z_{ij}^k \ \textit{для любых} (\alpha,\beta),\ k \in K, \qquad (4.14)$$

$$M_\alpha^k \ge \sum_{\{i,j \mid (\alpha,\beta)\in\Gamma_{ij}^k\}} z_{ij}^k \ \textit{для любых}\ \alpha,\ k \in K. \qquad (4.15)$$

Двойственной к задаче (4.1)-(4.3), (4.6), (4.14), (4.15) является задача

$$\sum_{i\in I} \Pi_i\left(p_i^1,\dots,p_i^m\right) + \sum_{j\in J} \pi_j\left(\hat{p}_j^1,\dots,\hat{p}_j^m\right) +$$

$$+\sum_{k\in K} \sum_{(\alpha,\beta)\in\Gamma_{ij}^k} t_{\alpha\beta}^k V_{\alpha\beta}^k + \sum_{(\alpha,\beta)\in\Gamma_{ij}^k} \lambda_\alpha^k M_\alpha^k \to \min \qquad , \ (4.16)$$

$$p_i^k \le \hat{p}_j^k + \sum_{(\alpha,\beta)\in\Gamma_{ij}^k} \left(t_{\alpha\beta}^k + \tilde{t}_{\alpha\beta}\right) + \sum_{(\alpha,\beta)\in\Gamma_{ij}^k} \left(\lambda_\alpha^k + \tilde{\lambda}_\alpha\right), \qquad (4.17)$$

$$p_i^k \ge 0, \hat{p}_j^k \ge 0,\ t_{\alpha\beta}^k \ge 0,$$

$$\lambda_\alpha^k \geq 0 \quad \left(i \in I, j \in J, k \in K, (\alpha,\beta) \in \Gamma_{ij}^k\right). \qquad (4.18)$$

Структура множителей Лагранжа в задаче (4.16)-(4.17) допускает перекрёстное субсидирование. Расчёт множителей может регулироваться за счёт выделения квот на перевозку различных товарных групп.

Подчеркнём, что перекрестное субсидирование позволяет не только регулировать уровень перевозок низко прибыльных товаров, но и позволяет, вообще говоря, выполнить социально значимые перевозки (4.10).

Отметим, что при решении задачи (4.1)-(4.3), (4.6), (4.14), (4.15) экстремум следует искать так же за счет выбора маршрутов перевозок $\Gamma_{ij}^k$.

# 5. Об одном подходе к анализу инвестиционного комплекса в системе железнодорожных грузоперевозок

Рассмотрим вопрос о распределении доходов от результатов экономической деятельности в конкурентном равновесии

$$\left\{X_i^k, Y_j^k, z_{ij}^k, p_i^k, \hat{p}_j^k \,\middle|\, i \in I, j \in J, k \in K\right\}$$

в модели железнодорожных грузоперевозок. Доход от деятельности рассматриваемой системы потребителей и производителей равен

$$\sum_{i \in I} F_i\left(X_i^1, \ldots, X_i^m\right) - \sum_{j \in J} G_j\left(Y_j^1, \ldots, Y_j^m\right).$$

Часть этого дохода

$$\sum_{i \in I, j \in J, k \in K} \left( \sum_{(\alpha,\beta) \in \Gamma_{ij}^k} \tilde{t}_{\alpha\beta} + \sum_{(\alpha,\beta) \in \Gamma_{ij}^k} \tilde{\lambda}_\alpha \right) z_{ij}^k \qquad (5.1)$$

идёт на оплату услуг по железнодорожным грузоперевозкам. В то же время совокупная прибыль потребителей равна

$$\sum_{i \in I} \Pi_i\left(p_i^1, \ldots, p_i^m\right), \qquad (5.2)$$

а совокупная прибыль производителей равна

$$\sum_{j\in J}\pi_j\left(\hat{p}_j^1,...,\hat{p}_j^m\right). \tag{5.3}$$

По теореме Фенхеля [12, с. 46-47] оптимальные значения функционалов в задаче (4.1)-(4.6) и двойственной к ней задаче (4.7)-(4.9) равны. Откуда следует, что коммуникационные ограничения пропускной способности (4.5), (4.6), влияющие на множители Лагранжа $\left\{t_{\alpha\beta},\lambda_\alpha\right\}$ и, тем самым, на разность цен в пунктах потребления и производства, порождают посредническую прибыль:

$$\sum_{(\alpha,\beta)\in\Gamma_{ij}^k}t_{\alpha\beta}V_{\alpha\beta}+\sum_{(\alpha,\beta)\in\Gamma_{ij}^k}\lambda_\alpha M_\alpha=\sum_{i\in I}F_i\left(X_i^1,...,X_i^m\right)-$$

$$-\sum_{j\in J}G_j\left(Y_j^1,...,Y_j^m\right)-\sum_{i\in I,\,j\in J,\,k\in K}\left(\sum_{(\alpha,\beta)\in\Gamma_{ij}^k}\tilde{t}_{\alpha\beta}+\sum_{(\alpha,\beta)\in\Gamma_{ij}^k}\tilde{\lambda}_\alpha\right)z_{ij}^k-$$

$$-\sum_{i\in I}\Pi_i\left(p_i^1,...,p_i^m\right)-\sum_{j\in J}\pi_j\left(\hat{p}_j^1,...,\hat{p}_j^m\right).$$

Из условий 6 и 7 в определении конкурентного равновесия следует, что

$$\sum_{(\alpha,\beta)\in\Gamma_{ij}^k}t_{\alpha\beta}V_{\alpha\beta}+\sum_{(\alpha,\beta)\in\Gamma_{ij}^k}\lambda_\alpha M_\alpha=$$

$$=\sum_{i\in I,j\in J,k\in K}\left(\sum_{(\alpha,\beta)\in\Gamma_{ij}^k}t_{\alpha\beta}+\sum_{(\alpha,\beta)\in\Gamma_{ij}^k}\lambda_\alpha\right)z_{ij}^k. \tag{5.4}$$

Увеличивая тарифы на услуги по грузоперевозкам с системы тарифов $\left\{\tilde{t}_{\alpha\beta},\tilde{\lambda}_\alpha\right\}$, определяемой по стоимости издержек на обслуживание грузоперевозок, до уровня $\left\{\tilde{t}_{\alpha\beta}+t_{\alpha\beta},\tilde{\lambda}_\alpha+\lambda_\alpha\right\}$, ОАО «РЖД» могло бы получить финансовые ресурсы для реализации инвестиционных программ, не изменяя при этом материальных потоков $\left\{X_i^k,Y_j^k,z_{ij}^k\middle|i\in I,j\in J,k\in K\right\}$.

Отметим, что реализация инвестиционных проектов по увеличению пропускной способности коммуникационных ограничений $\{M_\alpha, V_{\alpha\beta}\}$ (сортировочных станций или железнодорожных путей) приводит к уменьшению доходов (5.4). Поэтому для реализации этого предложения потребуется разделить финансовые потоки (5.1) и (5.4) и создать самостоятельную организационную структуру, финансирующую инвестиционные проекты по развитию системы железнодорожного транспорта и получающую денежные средства от ОАО «РЖД» (5.4) и государства. Реализация инвестиционных проектов по увеличению пропускной способности транспортной системы приводит к уменьшению размеров денежного потока (5.4), но увеличивает оптимальное значение функционала в задаче (4.1)-(4.6), а значит, совокупные прибыли потребителей (5.2) и (5.3). Было бы логично, если бы часть увеличивающихся за счёт этого поступлений с налога на прибыль государство направляло бы на инвестиции в систему железнодорожного транспорта.

В экстремальной задаче (4.7)-(4.9) для анализа тарифов на грузоперевозки множители Лагранжа $\{t_{\alpha\beta}, \lambda_\alpha\}$ являются оценками экономического эффекта от инвестиций в увеличение пропускной способности железной дороги $V_{\alpha\beta}$ по перевозке грузов из узловой станции $\alpha$ на узловую станцию $\beta$ и пропускной способности $M_\alpha$ узловой станции $\alpha$ по приёму и отправлению грузов. Решение задач (4.1)-(4.6) при различных значениях $V_{\alpha\beta}$, $M_\alpha$ позволяет выделить субъектов рассматриваемой экономической системы, заинтересованных в увеличении пропускных способностей и реализации соответствующих инвестиционных проектов.

# 6. Эволюционная интерпретация конкурентного равновесия

В связи с изучением вычислительных аспектов конкурентного равновесия полезно будет привести прямой метод получения равновесия и равновесных цен, имеющий содержательную эволюционную интерпретацию [15]. Далее мы будем отчасти следовать работам [15,16].

Предварительно “раздуем” исходный граф, считая, что каждой вершине $\alpha$ (см. ограничения (4.5’)) исходного графа соответствует дополнительное ребро (дуга): все ребра, входящие в эту вершину, входят в начало этой дуги, а все выходящие из этой вершины ребра выходят из конца этой дуги. Тогда ограничения на пропускную способность будут только у ребер графа. В раздутом графе также необходимо искусственно ввести один источник, который нужно соединить дугами с пунктами производства, и один сток, который нужно соединить дугами с пунктами потребления (ограничений на пропускные способности этих дуг нет). Источник и сток соединим дополнительным ребром с бесконечно большой пропускной способностью и нулевыми затратами.

Будем далее рассматривать “раздутый” граф $G = \langle V, E \rangle$. Обозначим через $P$ все пути (маршруты) в этом графе из источника в сток. Положим $x_p^k$ - объем перевозимого (в единицу времени) товара типа $k$ по пути $p \in P$. Определим (не ограничивая общности) функции затрат (для $k$-го товара) при прохождении вспомогательных ребер $\{\tilde{e}\}$, инцидентных источнику или стоку:

$$-p_i^k\left(X_i^1,\dots,X_i^m\right) = -\frac{\partial F_i\left(X_i^1,\dots,X_i^m\right)}{\partial X_i^k},$$

$$\hat{p}_j^k\left(Y_j^1,\dots,Y_j^m\right)=\frac{\partial G_j\left(Y_j^1,\dots,Y_j^m\right)}{\partial Y_j^k}, \tag{6.1}$$

где $\left\{X_i^k\right\}$, $\left\{Y_j^k\right\}$ аффинно выражаются через $\left\{x_p^k\right\}$. Определим функции затрат всех остальных ребер, на которые есть ограничения пропускной способности

$$\tau_e^\mu\left(f_e\right)=\bar{t}_e\cdot\left(1+\gamma\cdot\left(\frac{f_e}{\bar{f}_e}\right)^{\frac{1}{\mu}}\right), \tag{6.2}$$

где $\gamma>0$, $\mu>0$ – произвольный постоянные; $f_e$ – поток товаров (всех типов) на дуге $e\in E$ (аффинно выражается через $\left\{z_{ij}^k\left(\tilde{\Gamma}_{ij}\right)\right\}$, которые, в свою очередь, аффинно выражаются через $\left\{x_p^k\right\}$); $\bar{t}_e=\tilde{t}_{\alpha\beta}$ или $\bar{t}_e=\tilde{\lambda}_\alpha$, аналогично $\bar{f}_e=V_{\alpha\beta}$ или $\bar{f}_e=M_\alpha$ в зависимости от типа ребра.

Рассмотрим (популяционную) игру (загрузок) [15]. Пусть имеется много агентов, готовых перевозить из источника в сток товары разных типов (считаем, что в источнике неограниченно много товаров). Ничего не перевозя, агент имеет затраты ноль (прибыль ноль). Игра повторяется во времени. До тех пор, пока агентам выгодно перевозить товары, они будут вовлекаться в процесс перевозок, выбирая наиболее выгодные маршруты и товары для перевозок. Для описания этого эволюционного процесса есть много различных подходов [15], приводящих к одному и тому же результату на больших временах. А именно: система сойдется к равновесию (иногда называемому равновесием Нэша-Вардропа [16]), которое (ввиду того, что рассматривается игра загрузок, следовательно (см. Розенталь, 1973; Мондерер-Шэпли, 1996), потенциальная игра [15,17]) будет определяться минимумом (по $\left\{x_p^k\right\}\geq 0$) потенциальной функции этой системы $\Psi\left(\left\{x_p^k\right\}\right)$. Эта функция – функция Ляпунова различных естественных эволюционных динамик [15]:

$$\sum_{j\in J} G_j\left(\left\{Y_j^k\left(\left\{x_p^k\right\}\right)\right\}\right)-\sum_{i\in I} F_i\left(\left\{X_i^k\left(\left\{x_p^k\right\}\right)\right\}\right)+\sum_{e\in E\setminus\{\tilde{e}\}}\int_0^{f_e\left(\left\{x_p^k\right\}\right)}\tau_e(z)\,dz, \quad (6.3)$$

в нашем случае является выпуклой.

Не сложно проверить, что если $\mu \to 0+$, то (независимо от выбора $\gamma > 0$) задача (6.3) в точности переходит в задачу (4.1') – (4.6'). Причем

$$\tau_e^{\mu}\left(f_e\left(\left\{x_p^k(\mu)\right\}\right)\right)\xrightarrow[\mu\to 0+]{}\begin{cases}\tilde{t}_{\alpha\beta}+t_{\alpha\beta}\\ \tilde{\lambda}_{\alpha}+\lambda_{\alpha}\end{cases}$$

в зависимости от типа ребра $e \in E$. Здесь $t_{\alpha\beta}$ и $\lambda_\alpha$ – те же самые (двойственные множители), что и в теореме 4.2. С учетом формулы (6.1) имеем полный набор прямых переменных и всех двойственных множителей в теореме 4.2.

Можно показать [16], что от выбора монотонно возрастающей гладкой функции загрузок ребер вида (6.2) приведенный результат зависеть не будет. Единственное, что дополнительно нужно предполагать относительно функций $\tau_e^{\mu}(f_e)$, это

$$\tau_e^{\mu}(f_e)\xrightarrow[\mu\to 0+]{}\begin{cases}\bar{t}_e, & 0\le f_e<\bar{f}_e\\ \left[\bar{t}_e,\infty\right), & f_e=\bar{f}_e\end{cases},$$

Например, можно взять

$$\tau_e^{\mu}(f_e)=\bar{t}_e\cdot\left(1-\mu\ln\left(1-f_e/\bar{f}_e\right)\right).$$

Имитационная логит динамика (равномерно по $\mu \ge 0$) из [15] приводит к равновесию, т.е. решению задачи $\Psi\left(\left\{x_p^k\right\}\right)\to\min\limits_{\left\{x_p^k\right\}\ge 0}$.

Полезно также отметить, что если агенты ограниченно-рациональны [18], то, например, при обычной логит динамике (с параметром $\eta > 0$) [15] система сойдется не к равновесию, которое будет с точки зрения каждого агента равновесием Нэша, а к так называемому стохастическому равновесию. Это стохастическое равновесие, по сути, есть инвариантная мера некоторого эргодического марковского процесса, описывающего логит динамику. С ростом числа агентов $N$ (объема перевозимых товаров) и (или) при $\eta \to 0+$ эта инвариантная мера

$$\sim \exp\left( -N \cdot \left( \Psi\left(\left\{x_p^k\right\}\right) + \eta \sum_{k \in W} \sum_{p \in P} x_p^k \ln x_p^k \right) \Big/ \eta \right)$$

будет экспоненциально концентрироваться в малой окрестности решения задачи выпуклой оптимизации:

$$\Psi\left(\left\{x_p^k\right\}\right) + \eta \sum_{k \in W} \sum_{p \in P} x_p^k \ln\left(x_p^k\right) \to \min_{\left\{x_p^k\right\} \geq 0}, \tag{6.4}$$

которая является $\eta$-регуляризацией задачи (6.3) (в смысле расстояния Кульбака-Лейблера) или “энтропийной регуляризацией”. Здесь мы считаем, не ограничивая общности, что $\sum_{k \in W} \sum_{p \in P} x_p^k = 1$, интерпретируя $\left\{x_p^k\right\}$ как соответствующие доли. Важно отметить, что здесь функция $\Psi$ будет отличаться на шкалирующий множитель от функции $\Psi$, определенной в формуле (7.3).

Можно показать, что если $\eta \to 0+$ (при уже состоявшемся предельном переходе $N \to \infty$), то (всегда единственное в виду сильной выпуклости функционала) решение задачи (6.4) сходится к решению задачи (6.3) в случае единственности последнего. Если же решение задачи (6.3) не единственно, то отмеченный предельный переход представляет собой содержательно интерпретируемый способ отбора единственного равновесия, которое с большой вероятностью реализуется на практике [16].

Для численного решения задачи (6.3) можно использовать метод Франка–Вульфа (стандартный подход для такого класса задач), но при больших размерах графа более эффективным может оказаться рандомизированный метод из [16]. Оба эти метода могут по ходу итераций параллельно формировать решение двойственной задачи (двойственные множители), и оба могут быть проинтерпретированы содержательным образом. При этом из известных нам способов упомянутые способы самые эффективные для рассматриваемого класса задач.

В заключение отметим, что при фиксированных $\left\{X_i^k\right\}$, $\left\{Y_j^k\right\}$ к потенциальной игре загрузок, приводящей к задаче (6.3), эффективно применим механизм Викри–Кларка–Гроуса [17] (VCG mechanism), который говорит о том, как нужно выстраивать ценовую (тарифную) политику (определять дополнительные наценки $t_{\alpha\beta}$ и $\lambda_\alpha$), чтобы возникающее конкурентное равновесие в новой игре соответствовало “системному оптимуму” в исходной постановке. Применительно к рассматриваемой нами ситуации это подробно сделано в работах [19,20]. Интересно отметить, что при $\mu \to 0+$ выдаваемые этим механизмом тарифы $t_{\alpha\beta}$ и $\lambda_\alpha$ оказываются в точности такими, как в теореме 4.2.

## 7. Заключение

Проведённый на языке математических моделей анализ проблемы формирования тарифов на железнодорожные грузоперевозки позволяет сделать ряд практически полезных выводов.

1. Ограничения пропускной способности увеличивают разрыв между ценами в пунктах производства и ценами в пунктах потребления и порождают посредническую прибыль.
2. Ограничение на рост тарифов не позволяет системе железнодорожных грузоперевозок участвовать в извлечении посреднической прибыли и порождает недостаток финансовых средств на увеличения пропускной способности, модернизацию и обновление основных средств.
3. Привлечение частных инвестиций в обмен на привилегии в извлечении посреднической прибыли не является «системным решением проблемы», поскольку увеличение пропускной способности железнодорожной сети приводит к уменьшению разницы между ценами в пунктах производства и потребления и уменьшает посредническую прибыль.
4. Модельные примеры показывают, что максимизация посреднической прибыли может приводить к существенному сокращению объёмов производства и снижению эффективности грузовых потоков.
5. В отличие от агентов, заинтересованных в извлечении посреднической прибыли, общество в целом (государство), а также некоторые производители и потребители заинтересованы в увеличении пропускной способности грузоперевозок.

Вычислимые модели взаимодействия экономических агентов с учётом коммуникационных ограничений могут выявить инвесторов, заинтересованных в увеличении пропускной способности и модернизации железнодорожных грузоперевозок.

# 8. Литература

# Оглавление